\documentclass[12pt,reqno]{amsart}
\usepackage{amsmath, amscd}
\usepackage{amsfonts, amssymb}
\usepackage{mathtools}
\usepackage{xcolor}
\usepackage{mathrsfs}
\usepackage{geometry}
\usepackage{nicefrac}
\usepackage[all]{xy}

\usepackage[all]{xy}

\numberwithin{equation}{section}

\def\1#1{\overline{#1}}
\def\2#1{\widetilde{#1}}
\def\3#1{\widehat{#1}}
\def\4#1{\mathbb{#1}}
\def\5#1{\frak{#1}}
\def\6#1{{\mathcal{#1}}}

\newcommand{\R}{\mathbb R}

\newcommand{\C}{\mathbb C}

\newcommand{\B}{\mathbb B}

\newcommand{\Ha}{\mathbb H}

\def\v{\varphi}
\def\Re{{\sf Re}\,}

\newtheorem{theorem}{Theorem}[section]

\newtheorem{corollary}[theorem]{Corollary}
\newtheorem{question}[theorem]{Question}
\newtheorem{conjecture}[theorem]{Conjecture}
\theoremstyle{definition}

\theoremstyle{remark}
\newtheorem{remark}[theorem]{Remark}
\numberwithin{equation}{section}

\newcommand{\mcite}[1]{\csname b@#1\endcsname}

\long\def\REM#1{\relax}

\newcommand{\Maponto}
{\xrightarrow{\hbox{\lower.2ex\hbox{$\scriptstyle \smash{\mathsf{onto}}$}}\,}}
\newcommand{\Mapinto}
{\xrightarrow{\hbox{\lower.2ex\hbox{$\scriptstyle \smash{\mathsf{into}}$}}\,}}

\title[Abstract basins]{Abstract basins of attraction}
\author[L. Arosio]{Leandro Arosio$^\ast$}

\address{Dipartimento Di Matematica\\
Universit\`{a} di Roma \textquotedblleft Tor Vergata\textquotedblright\ \\
Via Della Ricerca Scientifica 1, 00133 \\
Roma, Italy} \email{arosio@mat.uniroma2.it}

\thanks{$^{*}$ Supported by the ERC grant ``HEVO - Holomorphic Evolution Equations'' n. 277691.}

\date\today

\begin{document}

\maketitle

\begin{abstract}
Abstract basins appear naturally in  different areas of several complex variables. In this survey we want to describe three different topics in which they play an important role, leading to interesting  open problems. 
\end{abstract}

\section{The construction of abstract basins}

In recent years  strong links among three different areas of research in several complex variables were discovered. It was indeed shown that 
the Bedford conjecture, the Loewner theory and  the theory of models for holomorphic self-maps revolve around a common concept introduced by Forn\ae ss--Stens\o nes \cite{Fo-St} in 2004: the abstract basin of attraction.
 Let $\B^q$ denote the open unit ball in $\C^q$, and let $(\varphi_{n}\colon \B^q\to \B^q)_{n\geq 0}$ be a family of univalent  (holomorphic injective) self-maps. We can think of this family as a dynamical system whose evolution law may change in time, and we call it a   {\sl  non-autonomous univalent dynamical system}. If $\varphi_n=\varphi$ for all $n\geq 0$, then we call it {\sl autonomous}.
 If $0\leq n\leq m$ we denote by $\varphi_{n,m}$ the composition $\varphi_{m-1}\circ \varphi_{m-2}\circ \cdots \circ \varphi_{n}$. 
 To construct the abstract basin we take the direct limit of  $(\varphi_{n})$, that is, we consider the following equivalence relation on the product  $\B^q\times \mathbb{N}$: let $0\leq n\leq m$, then $(x,n)\simeq (y,m)$ if and only if $\varphi_{n,m}(x)=y$. The set $\Omega=\B^q\times \mathbb{N}/_\sim$ is the {\sl abstract basin of attraction (or the tail space) of $(\varphi_{n})$}. It  comes naturally endowed with a family of mappings $(f_n\colon \B^q\to \Omega)_{n\geq 0}$ which satisfy $$f_n=f_m\circ \v_{n,m},\quad 0\leq n\leq m,$$ and which give a complex structure to $\Omega$. The abstract basin $\Omega$ satisfies the following universal property: if $\Lambda$ is a complex manifold, and $(g_n\colon \B^q\to \Lambda)_{n\geq 0}$ is a family of holomorphic mappings which satisfy 
 \begin{equation}\label{minnie}
 g_n=g_m\circ \v_{n,m},\quad 0\leq n\leq m,
 \end{equation}
  then there exists a holomorphic mapping $\Psi\colon \Omega\to \Lambda$ such that $f_n=\Psi\circ g_n$ for all $n\geq 0$. The mapping $\Psi$ is univalent if and only if  the mappings $g_n$ are univalent for all $n\geq 0$, and its image is the domain $\bigcup_{n\geq 0}g_n(\B^q)\subset \Lambda$.
\begin{remark}\label{evo}
The same construction works for a non-autonomous univalent dynamical system with $\R^+$ as an index set, that is, a  family of univalent mappings $(\varphi_{s,t}\colon \B^q\to\B^q)_{0\leq s\leq t}$ satisfying the {\sl evolution equation} 
\begin{equation}
\v_{u,t}\circ\v_{s,u}=\v_{s,t},\quad 0\leq s\leq u\leq t.
\end{equation}
If $t\mapsto \v_{s,t} (z)$ is locally lipschitz, uniformly on compacta in $z$, then $(\varphi_{s,t})$ is called an {\sl evolution family}. 
\end{remark}

We know by  construction that the basin $\Omega$ is the growing union of the subdomains $f_n(\B^q)$ which are biholomorphic to $\B^q$, however, the complex structure of  $\Omega$ may  be very complicated. Using a striking example due to Forn\ae ss \cite{Fornaess} it is easy to see that there exists an abstract basin which is not Stein and which is not biholomorphic to a domain of $\C^q$.

The name ``abstract basin of attraction'' deserves an explanation. Assume that there exists a non-autonomus dynamical system given by a family of automorphisms $(\Phi_n\colon X\to X)$ and a univalent mapping $h\colon \B^q\to X$ such that for all $n\geq 0$ we have $$h\circ \v_n=\Phi_n\circ h.$$ 
For all $n\geq0$, let $\Phi_{0,n}\coloneqq \Phi_{n-1}\circ \Phi_{n-2}\circ \cdots \circ \Phi_{0}$.
Then the family $g_n\coloneqq \Phi_{0,n}^{-1}\circ h$ satisfies equation (\ref{minnie}), and thus,
 by the universal property, the abstract  basin of attraction $\Omega$ of $(\v_n)$ is biholomorphic to the domain of $X$ given  by $$\bigcup_{n\geq 0}g_n(\B^q)=\{x\in X\colon \Phi_{0,n}(x)\in h(\B^q) \ \mbox{eventually}\}.$$
Assume now that the origin is attracting for the dynamical system  $(\v_n)$, that is $\v_{0,n}(z)\to 0$ for all $z\in \B^q$.  Then  $ \Phi_{0,n}(x)\to h(0) $ for all $x\in h(\B^q)$, which implies that the abstract basin of attraction
$\Omega$ of $(\v_n)$ is biholomorphic to  the ``actual''   basin of attraction of $(\Phi_n)$ at the point $h(0)$, that is the domain of $X$ defined by
$$\{x\in X\colon \Phi_{0,n}(x)\to h(0)\}.$$

\section{Bedford's Conjecture}

Let  $\v\colon \B^q\to \B^q$ be a univalent mapping fixing the origin such that the spectrum of the differential at the origin $d_0\v$ is contained in the punctured disc $\Delta\smallsetminus\{0\}$. Then the autonomous univalent dynamical system associated with $\varphi$ is attracting at the origin, and its abstract basin of attraction is biholomorphic to $\C^q$ by the Poincar\'e--Dulac theory (see e.g. \cite{RR}).
 
It is natural to ask whether this holds true if we consider non-autonomous dynamical systems $(\v_n\colon\B^q\to \B^q)$ whose contraction rate at the origin is uniformly bounded from above and from below. 
\begin{conjecture}[\cite{Fo-St}]\label{1}
Let $(\v_n\colon\B^q\to \B^q)$ be a non-autonomous univalent dynamical system such that $\v_n(0)=0$ for all $n\geq 0$ and that 
$$a\|z\|\leq \|\varphi_n(z)\|\leq b\|z\|,\quad n\geq 0, \ 0<a\leq b<1.$$
Then the abstract basin $\Omega$ of $(\v_n)$ is biholomorphic to $\C^q$.
\end{conjecture}

\begin{remark}
Condition $b<1$ ensures that $\v_{0,n}(z)\to 0$ for all $z\in \B^q$.
If we drop the condition $a>0$ there is the following counterexample by Forn\ae ss \cite{Fornaess-short}. Set $$\v_n\colon (z,w)\mapsto (z^2+a_n w, a_n z),$$ where $|a_0|<1$ and $|a_{n+1}|\leq a_n^2$. Then $\Omega$ admits a non-constant bounded plurisubharmonic function, and thus is not biholomorphic to $\C^2$. Such a domain is called a {\sl short $\C^2$}, since it shares several invariants with $\C^2$ without being biholomorphic to $\C^2$.
\end{remark}

Conjecture \ref{1}  is still open, and is deeply studied since Forn\ae ss--Stens\o nes \cite{Fo-St} proved  that it is stronger than   the well-known  {\sl Bedford conjecture}:
\begin{conjecture}[\cite{Bedford}]
Let $X$ be a complex manifold endowed with a Riemannian metric, and let $f\colon X\to X$ be an automorphism which acts hyperbolically on some invariant compact subset $K\subset X$. If $p\in K$, then  the stable manifold $\Sigma(p)$ is biholomorphic  to $\C^k$, where $k$ is the stable dimension.
\end{conjecture}
Jonsson--Varolin \cite{JV} showed that the Bedford conjecture is true for every $p$ in a subset of $K$ which is of full-measure with respect to any invariant probability measure on $K$, and Abate--Abbondandolo--Majer \cite{AAM} showed that if the Lyapounov exponents exist, then $\Sigma(p)$ is biholomorphic  to $\C^k$.

Forn\ae ss--Stens\o nes \cite{Fo-St} proved that $\Omega$ is always biholomorphic to a domain of $\C^q$. This implies in particular that $\Omega$ has to be Stein, and it is also easy to see that the Kobayashi pseudometric and pseudodistance of $\Omega$ are zero everywhere. The main approach to  Conjecture \ref{1} consists in adapting the Poincar\'e--Dulac method to the non-autonomous setting, and this is why
 arithmetic relations between $a$ and $b$ play a capital role. Wold \cite{Wold} showed that if $b^2<a$, then $\Omega$ is biholomorphic to $\C^q$. Abbondandolo--Majer    \cite{AM} showed that in $\B^2$ this condition can be loosen to $b^{29/14}<a$. Very recently  Peters--Smit \cite{Peters-Smit} obtained $b^{\frac{11}{5}}<a$, and, assuming all $d_0\v_n$ diagonal, $b^{3}<a$.
  
  Other interesting ``weak monotonicity'' relations among the eigenvalues of $d_0 \v_n$ are considered in \cite{Peters, AAM}. See \cite{everybody} for a survey on non-autonomous basins and the Bedford conjecture.

\section{Loewner's theory}
The Loewner PDE in the unit disc was introduced by Loewner \cite{Loewner} in 1923 while he was working on the Bieberbach conjecture. It was later developed by Kufarev \cite{Kuf1943} and Pommerenke \cite{Pommerenke} and  is today one of the principal tools in geometric function theory.
Recently Bracci--Contreras--D\'iaz-Madrigal   \cite{BCD1,BCD2} introduced a framework for a very general Loewner theory which works in complete hyperbolic manifolds (see also \cite{AB}). We consider the case of $\B^q$. They give the following definition, which is a generalization to the setting of continuous time of the concept of non-autonomous univalent dynamical system.
A {\sl Herglotz vector field} on $\B^q$ is a function $h(z,t)\colon \B^q\times \R^+\to \C^q$ such that 
\begin{enumerate}
\item $z\mapsto h(z,t)$ is a semicomplete holomorphic vector field,
\item $t\mapsto h(z,t)$ is measurable and locally bounded, uniformly on compacta in $z$.
\end{enumerate}

They prove that a Herglotz vector field $h$ is semicomplete in the sense that the ODE
$$\frac{dz(t)}{dt}=h(z,t)$$ has a solution flow given by an evolution family $(\v_{s,t}\colon \B^q\to \B^q)_{0\leq s\leq t}$. 

The Loewner PDE in several complex variables was studied by Pfaltzgraff, Graham, Duren,  Kohr, Hamada, and others (see \cite{Pf74,Pf75,DGHK,GHK1,GHK2}).
In \cite{ABHK}, the Loewner PDE was generalized to the setting of Herglotz vector fields in the following way:
\begin{equation}\label{PDE}
\frac{\partial f_t(z)}{\partial t}=-d_zf_t \ h(z,t), \quad {\mbox a.e.}\ t\geq 0, z\in \B^q,
\end{equation}
where the unknown $(f_t\colon \B^q\to \C^q)$ is a family of univalent mappings such that $t\mapsto f_t(z)$ is locally lipschitz, uniformly on compacta in $z$.

The following result shows that solutions exist and  are essentially unique  if we do not restrict ourselves to solutions with values in $\C^q$. 

\begin{theorem}[\cite{ABHK}]
The Loewner PDE (\ref{PDE}) admits a univalent solution $(f_t\colon \B^q\to \Omega)$, where $\Omega$ is
the abstract basin of attraction of the evolution family $(\v_{s,t})$ associated with the Herglotz vector field $h(z,t)$. If $(g_t\colon \B^q\to Q)$ is another solution, where $Q$ is a $q$-dimensional complex manifold, then there exists a holomorphic mapping $\Psi\colon \Omega\to Q$ such that $$g_t=\Psi\circ f_t.$$
\end{theorem}
The abstract basin $\Omega$ is called the {\sl Loewner range} of $h(z,t)$ (or of $(\v_{s,t})$).

\begin{remark}
This result transforms the analytic problem of finding a univalent solution for the Loewner PDE (\ref{PDE}) to the geometric problem of understanding whether the Loewner range $\Omega$ of $(\v_{s,t})$ is biholomorphic to a domain of $\C^q$. For example, in one variable, we know that $\Omega$ is non-compact and simply connected (since it is the growing union of discs). Thus, by the uniformization theorem, it has to be biholomorphic to $\C$ or to the disc $\Delta$. In either case it is biholomorphic to a domain of $\C$, and hence we obtain that  the Loewner PDE in one variable always admits a univalent solution, as proved by Contreras--D\'iaz-Madrigal--Gumenyuk \cite{one} with a different method.
\end{remark}

This embedding problem was solved in \cite{ABW} using two major tools: a result of Docquier--Grauert \cite{DocquierGrauert} which implies that for all $0\leq s\leq t$ 
the pair $(f_s(\B^q),f_t(\B^q))$ is Runge, and Anders\'en--Lempert theory.
\begin{theorem}
The Loewner range $\Omega$ of $(\v_{s,t})$ is biholomorphic to a domain of $\C^q$. As a consequence,  the Loewner PDE (\ref{PDE}) always admits a univalent solution $(f_t\colon \B^q\to \C^q).$
\end{theorem}

An interesting open question is to find conditions for a Herglotz vector field $h(z,t)$ which ensure that its Loewner range $\Omega$ is biholomorphic to $\C^q$. Some conditions are given in \cite{Ar1,Ar2,Ar3}.  We can also formulate a  continuous-time analogue of the Bedford conjecture. Recall that if $A$ is a linear endomorphism of $\C^q$, we denote  $m(A)\coloneqq \min\{\Re\langle Az,z\rangle\colon |z|=1\}$ and $k(A)\coloneqq \max\{\Re\langle Az,z\rangle\colon |z|=1\}$.

\begin{conjecture}\label{domanda}
Let $h(z,t)$ a Herglotz vector field on $\B^q$ of the form $h(z,t)=A(t)z+O(|z|^2)$. Assume that 
\begin{enumerate}
\item $m(A(t))>0$ for all $t\geq 0$ and $\int_0^\infty m(A(t))d(t)=\infty,$
\item there exists $\ell\in\R^+$ such that $\ell m(A(t))\geq k(A(t))$, for all $t\geq 0.$
\end{enumerate}
Then the Loewner range of $h(z,t)$ is biholomorphic to $\C^q$.

\end{conjecture}

Consider the evolution family $(\v_{s,t})_{0\leq s\leq t}$ associated with $h(z,t)$. Notice that any discretization of  the index set $\R^+$ gives as a result  a (discrete) non-autonomous dynamical system $(\v_{n,m})_{0\leq n\leq m}$. The abstract basins of the two families are easily seen to be biholomorphic. 
The assumptions in Conjecture \ref{domanda} allow to discretize the index set $\R^+$ in such a way that  $(\v_{n,m})_{0\leq n\leq m}$ satisfies the assumption of the Bedford conjecture. Thus Conjecture \ref{domanda} is weaker than the Bedford conjecture.

For a survey on  Loewner theory, see \cite{ABCD}.

\section{Models}

The idea of using representation models to understand the local dynamics of holomorphic self-maps goes back to the birth of complex dynamics itself, that is the introduction  in 1870 of the Schr\"oder equation \cite{S1,S2}.
Let $f\colon \Delta\to \Delta$ be a holomorphic mapping fixing  the origin.  Assume that $0<|f'(0)|<1$. Then the origin is an attracting fixed point, and the Schr\"oder equation is the following:
\begin{equation}
\sigma\circ f=f'(0)\circ \sigma,
\end{equation}
where  $\sigma\colon \Delta\to \C$ is an unknown holomorphic function. This equation was solved in $1884$ by K\"onigs \cite{Ko}, which showed that there  exists a  solution $\sigma$, which is unique if we impose $\sigma(0)=0,\ \sigma'(0)=1$.

If $f$ has no interior fixed points, then by the Denjoy--Wolff theorem there exists a point $a\in \partial \Delta$ called the {\sl Denjoy--Wolff point} such that $f^n$ converges to $a$ uniformly on compacta. Moreover we can define the {\sl dilation of $f$} as the following non-tangential limit:
$$\lim_{z\to a}f'(z)=\lambda\in (0,1].$$ 
The mapping $f$ is called  {\sl hyperbolic} iff $\lambda\in (0,1)$, and is called {\sl parabolic} iff $\lambda=1$. 

Let $\Ha$ denote the upper half-plane.
If $f$ is  hyperbolic, Valiron \cite{Va} proved in 1931 that there exists  $\sigma\colon \Delta\to\mathbb{H}$ such that $$\sigma\circ f=\frac{1}{\lambda}\sigma,$$ and any other solution is a positive multiple of $\sigma$. 
Notice that the growing union $\cup_{m\in \mathbb{N}}\, \lambda^m\sigma(\Delta)$ fills the whole half-plane $\mathbb{H}$.
If $f$ is   parabolic, Pommerenke--Baker \cite{Po1,BaPo} 
proved in 1979 that there exists $\sigma\colon \Delta\to\mathbb{C}$ such that 
$$\sigma\circ f=\sigma+1.$$
In this case we have two cases for the complex structure of the growing union  $\cup_{m\in \mathbb{N}}\, (\sigma(\Delta)-m)$.
Recall that, if $z_m\coloneqq f^m(z_0)$ is an orbit, its {\sl step} $s(f,z_0)$ is defined as $\lim_{m\to\infty} k_{\Delta}(z_m,z_{m+1})$, where $k_{\Delta}$ denotes the Poincar\'e distance of the disc $\Delta$. Such a limit exists thanks to the non-expansiveness of the Poincar\'e distance.  
We have the following dichotomy:
\begin{enumerate}
\item for any orbit $(z_m)$ we have $ s(f,z_0)= 0$ ({\sl zero-step}),
\item  for any orbit $(z_m)$ we have $ s(f,z_0)>0$ ({\sl nonzero-step}).
\end{enumerate}
The union   $\cup_{m\in \mathbb{N}}\, (\sigma(\Delta)-m)$  fills the whole $\mathbb{C}$ in the zero-step case, and is biholomorphic to $\Delta$ in the nonzero-step case.

Let now $f\colon \mathbb{B}^q\to \mathbb{B}^q$.  If the origin is an attracting  fixed point, then the Poincar\'e--Dulac theory applies and  one can solve a generalized Schr\"oder equation in several complex variables (see e.g. Rosay--Rudin \cite{RR}).
If there are no fixed points in $\mathbb{B}^q$,  then as in one variable there exists a {\sl Denjoy--Wolff point}  $a\in \partial \mathbb{B}^q$,  such that $f^n\to a$ uniformly on compact subsets. We can define the {\sl dilation of $f$}  as $$\liminf_{z\to a}\frac{1-\|f(z)\|}{1-\|z\|}=\lambda\in (0,1].$$
Again the mapping $f$ is called {\sl hyperbolic} iff $\lambda\in (0,1)$, and is called {\sl parabolic} iff $\lambda=1$. 
Zero-step and nonzero-step are defined as in the disc using the Kobayshi distance instead of the Poincar\'e distance (notice however that it is not a dichotomy anymore).

There are several generalizations by Bracci, Poggi-Corradini, Gentili, Bayart, Jury (see \cite{BG,BGP,jury,bayart}) of the Valiron and Abel equations in the unit ball which require additional regularity at the Denjoy--Wolff point $a\in \mathbb{B}^q$. All are obtained by scaling limit arguments.
In 1981 Cowen \cite{Cowen} unified the Schr\"oder, Valiron and Abel equations in a unique framework, introducing the concept of {\sl model} (without naming it). 
Models in several complex variables were recently introduced in \cite{ABmod}. 
If $f\colon \B^q\to \B^q$ is a univalent mapping, a {\sl semi-model for $f$} is given by  a complex manifold  $\Omega$ ({\sl the base space}), an automorphism $\psi$ of $\Omega$ and  a holomorphic mapping $h\colon \mathbb{B}^q\to \Omega$ such that
$h\circ f=\psi\circ h$.
We also assume   $\Omega=\cup_{m\in \mathbb{N}}\, \psi^{-m}(h(\mathbb{B}^q)).$
If $h$ is univalent, then we call $(\Omega,\psi,h)$  a {\sl model for $f$}.

\begin{theorem}[\cite{ABmod}]
Every univalent mapping $f\colon \B^q\to \B^q$ admits a model $(\Omega, h, \psi)$, where  $\Omega$ is the abstract basin of attraction of the autonomous univalent dynamical system associated with $f$. Moreover, if $(\Lambda, \ell, \phi)$ is any semi-model for $f$, then there exists a surjective holomorphic map $\Psi: \Omega\to \Lambda$ 
such that  the following diagram commutes:
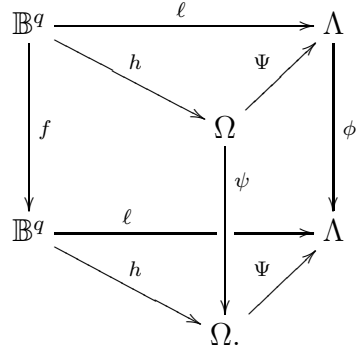
\SelectTips{xy}{12}
\[ \xymatrix{\B^q \ar[rrr]^\ell\ar[rrd]^h\ar[dd]^f &&& \Lambda \ar[dd]^\phi\\
&& \Omega \ar[ru]^\Psi \ar[dd]^(.25)\psi\\
\B^q\ar'[rr]^\ell[rrr] \ar[rrd]^h &&& \Lambda\\
&& \Omega \ar[ru]^\Psi.}
\]
\end{theorem}

The complex structure of $\Omega$ is not known in general. However, we can single out a semi-model on a possibly lower dimensional ball $\B^k$ which contains all the Kobayashi pseudodistance information of the model. 
Since $\Omega$ is the growing union of domains which are biholomorphic to $\B^q$, there exists, by a result of  Forn\ae ss--Sibony \cite{FS},  a surjective holomorphic submersion $r\colon \Omega\to \B^k$, where $0\leq k\leq q$, on whose fibers the Kobayashi pseudo-distance of $\Omega$ vanishes. 
This implies that the automorphism $\psi$ preserves the fibers, inducing an automorphism $\tau$ of $\B^k$ such that  the following diagram commutes:
$$
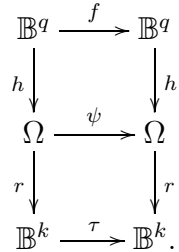
\xymatrix{\B^q\ar[r]^{f}\ar[d]_{h}& \B^q\ar[d]^{h}\\
\Omega\ar[r]^{\psi}\ar[d]_{r}& \Omega\ar[d]^{r}\\
\B^k\ar[r]^{\tau}&\B^k.}$$

Set $\ell\coloneqq r\circ h$.
The semi-model $(\B^k,\ell,\tau)$ is called the {\sl canonical semi-model}: every other semi-model with Kobayashi hyperbolic base space is a factor of this one.

Let us discuss what happens in the unit disc. If $f$ is hyperbolic with dilation $\lambda$, then by the Valiron equation it has a model with base space $\Delta$, and the automorphism $\psi$ is  hyperbolic  with dilation $\lambda$. The canonical semi-model coincides with the model.
If $f$ is parabolic nonzero-step, then by the Abel equation it has a model with base space $\Delta$, and the automorphism $\psi$ is  parabolic. The canonical semi-model coincides with the model.
If $f$ is parabolic zero-step the base space of the model is $\C$ and thus the base space of the canonical semi-model is a single point.
In \cite{ABmod}, the first two cases are generalized to several complex variables. 
\begin{theorem}
Let $f\colon \B^q\to \B^q$ be a hyperbolic (resp. parabolic nonzero-step) univalent mapping with dilation $\lambda$. Let  $(\B^k,\ell,\tau)$ be its canonical semi-model. Then $k\geq 1$ and $\tau$ is an hyperbolic (resp. parabolic) automorphism with the same dilation $\lambda$.
\end{theorem}

\begin{corollary}
Let $f\colon \B^q\to \B^q$ be a hyperbolic univalent mapping with dilation $\lambda$.  Then there exists a
holomorphic  solution  $\Theta\colon \B^q\to\Ha$ to the Valiron equation
$$\Theta\circ f=\frac{1}{\lambda}\Theta.$$
\end{corollary}

The following questions are open.
\begin{question}
 It would be natural to conjecture the following dichotomy (which is true in one variable and for linear fractional mappings of $\B^q$): if $f\colon \B^q\to \B^q$ is  a  parabolic univalent mapping and $(\B^k,\ell,\tau)$ is its canonical semi-model, then either $\tau$ is a parabolic automorphism, or $k=0$ and the canonical semi-model is trivial.

It is known  that the automorphism $\tau$ cannot be hyperbolic, however it is an open question whether a parabolic univalent mapping $f\colon \B^q\to \B^q$ can have a canonical semi-model $(\B^k,\ell,\tau)$ with $k>0$ and $\tau$ with an interior fixed point.
\end{question}
\begin{question}\label{q2}
If $f\colon \B^q\to \B^q$ is univalent, are the fibers of the holomorphic submersion $\Omega\to \B^k$ biholomorphic to $\C^{q-k}$? Is it true that $\Omega$ is biholomorphic to $\B^k\times \C^{q-k}$?
\end{question}

\begin{remark}
By a result of Forn\ae ss--Sibony \cite{FS}, if $k=q-1$, then the answer to Question \ref{q2} is affirmative.
\end{remark}

\end{document}